\documentclass[a4paper]{article}
\usepackage[left=4cm, right=4cm, bottom=3.5cm, top=3.5cm]{geometry}
\usepackage{graphicx}

\providecommand{\keywords}[1]{\textit{Keywords: } #1}

\usepackage[hyperfootnotes=false]{hyperref}
\usepackage{amsfonts, amsbsy, amssymb, amsmath}

\begin{document}

\def \cK {\mathcal{K}}

\title{Using Lagrangian descriptors to uncover invariant structures in Chesnavich's Isokinetic Model 
with application to roaming}

\author{Vladim{\'i}r Kraj{\v{n}}{\'a}k\footnote{School of Mathematics, University of Bristol, Bristol BS8 1TW, United Kingdom}, %
Gregory S. Ezra\footnote{Department of Chemistry and Chemical Biology, Baker Laboratory, Cornell University, Ithaca, NY 14853, USA}, %
Stephen Wiggins\footnotemark[1]}

\maketitle

\begin{abstract}
Complementary to existing applications of Lagrangian descriptors as an exploratory method, 
we use Lagrangian descriptors to find invariant manifolds in a system where some invariant 
structures have already been identified. In this case we use the parametrisation of a 
periodic orbit to construct a Lagrangian descriptor that will be locally minimised on its invariant manifolds. 
The procedure is applicable (but not limited) to systems with highly unstable periodic orbits, 
such as the isokinetic Chesnavich CH$_4^+$ model subject to a Hamiltonian isokinetic theromostat. Aside from 
its low computational requirements, the method enables us to study the invariant 
structures responsible for roaming in the isokinetic Chesnavich CH$_4^+$ model.

\end{abstract}

\section{Introduction}
\label{sec:intro}

This paper is concerned with the computation of stable and unstable manifolds of 
unstable periodic orbits using the method of Lagrangian descriptors (LDs) 
\cite{madrid2009ld, craven2016deconstructing, craven2017lagrangian, craven2015lagrangian, junginger2016transition, junginger2017chemical, feldmaier2017obtaining, junginger2016transition, feldmaier2017obtaining, patra2018detecting}.  
We illustrate this method by applying it to understand the phase space 
structure governing the roaming reaction mechanism 
\cite{suits2008, bowman2011roaming,Bowman2011Suits, BowmanRoaming, mauguiere2017roaming} in a setting where the system of interest 
is subjected to a nonholonomic constraint that enforces constant kinetic energy 
\cite{Dettmann96,Morriss98,Litniewski93,Morishita03}. 
First we give some background on this problem.

The roaming mechanism for chemical reactions was discovered in efforts to explain 
certain experimental data describing the photodissociation of formaldehyde \cite{zee:1664, townsend2004roaming, bowman}. 
Roaming provides a route to molecular dissociation products H$_2$ and CO that does not
involve passage through a conventional transition state.  
After excitation of the formaldehyde molecule (H$_2$CO) by a laser pulse, a single CH bond begins to stretch. 
Rather
than proceed directly to dissociation, one hydrogen atom rotates around the molecular 
fragment in a flat region of the potential energy surface; at a later stage this roaming hydrogen atom 
encounters the bound hydrogen atom and undergoes an abstraction reaction. 
The resulting H$_2$ molecule then separates from the CO fragment. This reaction is then said 
to occur by the roaming mechanism.

Since the pioneering formaldehyde studies roaming has been observed to 
occur in a number of chemical reactions, which are discussed in a number of  review articles 
\cite{suits2008, bowman2011roaming,Bowman2011Suits, BowmanRoaming, mauguiere2017roaming}.

We note that essentially all of these roaming reactions have been studied at constant total energy.  

Chesnavich developed an empirical model \cite{Chesnavich1986}  for 
the ion-molecule reaction CH$_4^+$ $\rightleftharpoons$ CH$_3^+$ $+$ H that contains 
the essential features of the roaming  mechanism. A detailed derivation of this model
can be found in \cite{ezra2019chesnavich}. Chesnavich's model describes the situation 
where a hydrogen atom separates from a rigid CH$_3^+$ core and, instead of dissociating, 
roams in a region of nearly constant potential only to return to the core.  
While Chesnavich's model does not accurately describe the intramolecular abstraction and subsequent dissociation, 
it has nevertheless provided significant insight into the roaming process, see, 
for example, \cite{mauguiere2014multiple, mauguiere2014roaming, krajnak2018phase, krajnak2018influence}.

The goal of this article is to extend  our previous study of the analogous roaming reaction mechanisms 
at constant kinetic energy  \cite{krajnak2019isokinetic}. 
We do this by applying a Hamiltonian isokinetic thermostat \cite{Dettmann96,Morriss98,Litniewski93,Morishita03} 
to Chesnavich's model. Use of the isokinetic thermostat means that we are effectively investigating 
roaming at constant temperature.
A  detailed study of Chesnavich's model subjected to the Hamiltonian isokinetic 
thermostat is given in \cite{krajnak2019isokinetic}. 
In that work it was found that certain of the periodic orbits that govern the roaming dynamics are 
highly unstable: 
their Lyapunov exponents are $\sim 10^{21}$. In such a situation detection of the periodic orbits, 
as well as the computation of their stable and unstable manifolds, becomes problematic using traditional approaches.

In this paper we show that the method of Lagrangian descriptors 
\cite{madrid2009ld, craven2016deconstructing, craven2017lagrangian, craven2015lagrangian, junginger2016transition, junginger2017chemical, feldmaier2017obtaining, junginger2016transition, feldmaier2017obtaining, patra2018detecting}
can be used to compute stable and unstable manifolds of such highly unstable periodic orbits. 
We applying this approach to Chesnavich's model subjected to a Hamiltonian isokinetic 
thermostat and show that the method of Lagrangian descriptors reveals the 
dynamical origin of the roaming mechanism in this setting.

This paper is outlined as follows: In Sec.\ \ref{sec:CMR} we introduce Chesnavich's Hamiltonian model  
for ion-molecule reaction and discuss the dynamical mechanism underlying roaming in 
terms of families of unstable periodic orbits and their associated invariant manifolds.
In Sec.\ \ref{sec: HIK} we discuss the Hamiltonian formulation of the isokinetic thermostat 
for Chesnavich's CH$_4^{+}$ model. In  Sec.\ \ref{sec:LDIM} we introduce 
the method of Lagrangian descriptors and use prior knowledge of periodic orbits to 
construct a Lagrangian descriptor for detecting stable and unstable invariant manifolds. 
We then apply the method to Chesnavich's model subjected to a Hamiltonian isokinetic 
thermostat to reveal the roaming mechanism in this system. Sec.\ \ref{sec:conc} concludes.

\section{Chesnavich's Model and Roaming}
\label{sec:CMR}
\subsection{Chesnavich's Model Hamiltonian}

The  CH$_4^+$ model due to Chesnavich  is a 2 degree of freedom Hamiltonian system  
comprised of a rigid CH$_3^+$ molecule (core) and a mobile H atom \cite{Chesnavich1986}. 
The system Hamiltonian is \cite{ezra2019chesnavich}
\begin{equation}
H(r,\theta,p_r, p_\theta) = \frac{1}{2} \frac{p_r^2}{\mu} + \frac{1}{2}p_\theta^2 \left(\frac{1}{\mu r^2}+\frac{1}{I_{CH_3}}\right) + U(r,\theta),
\label{eq:chesHam}
\end{equation}
where $(r, \theta, \phi)$ are polar coordinates describing the position of the 
H-atom in a body-fixed frame attached to the CH$_3^+$ core
(the coordinate $\phi$ is ignorable in this model).  
The reduced mass of the system is given by the expression 
$\mu=\frac{m_{CH_3}m_{H}}{m_{CH_3}+m_{H}}$, where $m_{H}=1.007825$ u and $m_{CH_3}=3m_{H}+12.0$ u, 
and the moment of inertia of the rigid body CH$_3^+$ has the value  $I_{CH_3}=2.373409$~u\AA$^2$.

The potential energy function $U(r,\theta)$ is made up  of a radial long range potential
energy term $U_{CH}$ and a short range potential $U_{coup}$ that models the short range anisotropy of the rigid CH$_3^+$ core:
\begin{equation}\label{eq:U}
 U(r,\theta ) = U_{CH} (r) + U_{coup} (r,\theta).
\end{equation}

\begin{figure}[ht]
 \centering
 \includegraphics[width=.5\textwidth]{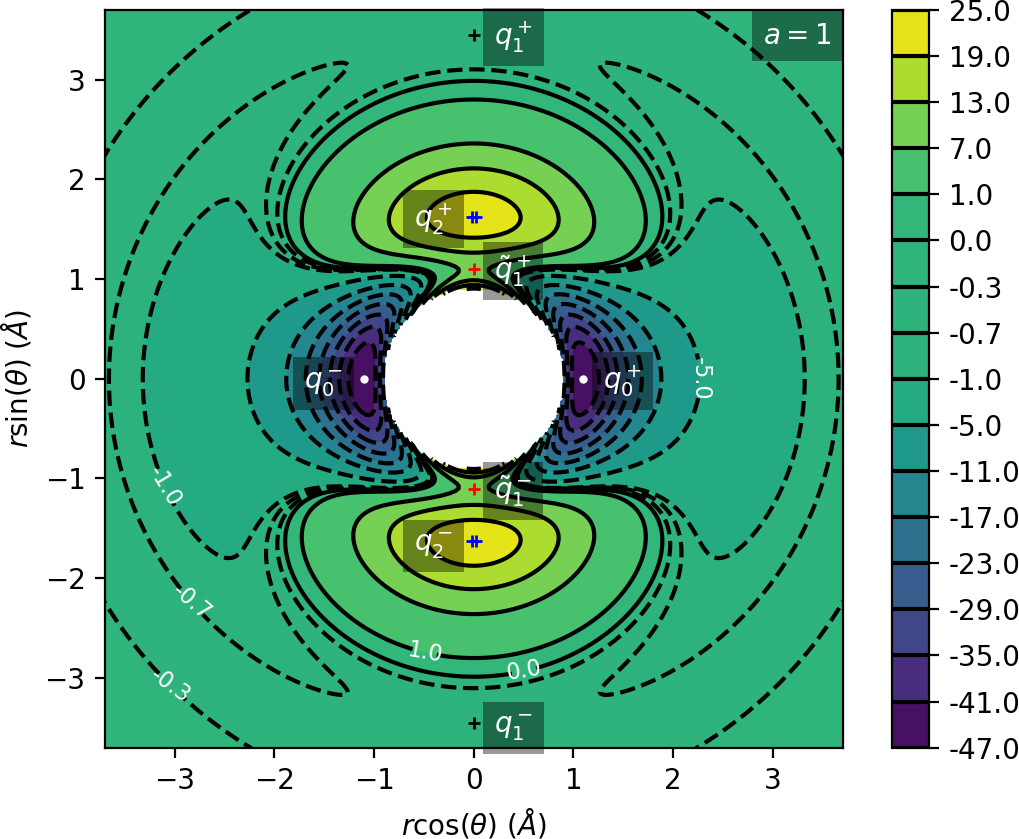}
 \caption{Contour plot of Chesnavich's potential energy surface $U$ for $a=1$. Dashed lines correspond to $U<0$, solid lines correspond to $U\geq0$. Contours correspond to values of potential shown on the colorbar to the right, with some values indicated in the plot. This figure is from \cite{krajnak2018influence}.
 }
 \label{fig:pot}
\end{figure}

The topography of the potential energy surface  is characterised by two deep wells that 
correspond to the  bound CH$_4^+$, two areas of high potential and a flat area to the outside 
of these features as shown in Fig. \ref{fig:pot}.

The long range potential has the form:
\begin{equation}\label{eq:UCH}
 U_{CH} (r) =  \frac{D_e}{c_1 - 6} \left( 2 (3-c_2) e^{c_1 (1-x)}  - \left( 4 c_2 - c_1 c_2 + c_1 \right) x^{-6} - (c_1 - 6) c_2 x^{-4} \right), 
\end{equation}
where $x = \frac{r}{r_e}$ and we take the  parameter values as  in the original work \cite{Chesnavich1986}.
The short range hindered rotor  potential $U_{coup}$ has the form:
\begin{equation}\label{eq:Ucoup}
 U_{coup} (r,\theta) = \frac{U_e e^{-a(r-r_e)^2}}{2} (1 - \cos 2 \theta ),
\end{equation}
where $U_e$ is the equilibrium barrier height. 
The distance at which the  transition occurs from rotation to vibration is 
determined by the parameter $a$ (in \AA$^{-2}$). 
Various values of $a$ have been considered in previous works. In particular,  
$a=1$ \cite{Chesnavich1986,mauguiere2014roaming,mauguiere2014multiple,krajnak2018phase}, 
$a=4$ \cite{Chesnavich1986,mauguiere2014roaming} and a range of values $0.7\leq a\leq 8$. \cite{krajnak2018influence}

The CH$_3^+$ core is a symmetric top in Chesnavich's model. 
Although the range of the coordinate $\theta$ is $0 \leq \theta \leq \pi$, 
in the planar (zero overall angular momentum) version of the model the range of $\theta$ 
is extended to $0 \leq \theta \leq 2 \pi$, and the potential has a four-fold symmetry:
\begin{equation}
U(r,\theta)=U(r,-\theta)=U(r,\pi-\theta)=U(r,\pi+\theta).
\label{eq:sym}
\end{equation}

The potential admits four pairs of equilibrium points pairwise related by symmetry \eqref{eq:sym}, 
as listed in Tab. \ref{table:equil} and shown in Fig. \ref{fig:pot}.
  
  \begin{table}[ht]
    \begin{center}
    \begin{tabular}{c|c|c|c|c}
    Energy (kcal mol$^{-1}$) & $r$ (\AA) & $\theta$ (radians) & Significance & Label \\
    \hline
    $-47$ & $1.1$ & $0$ & potential well & $q_0^+$ \\
    $-0.63$ & $3.45$ & $\pi/2$ & isomerisation saddle & $q_1^+$ \\
    $8$ & $1.1$ & $\pi/2$ & isomerisation saddle & $\widetilde{q}_1^+$ \\
    $22.27$ & $1.63$ & $\pi/2$ & local maximum & $q_2^+$ \\
    \end{tabular}
    \end{center}
    \caption{Equilibrium points of the potential $U(r, \theta)$.}
    \label{table:equil}
   \end{table}

\subsection{Roaming in Chesnavich's Model}
\label{sec:roamingChes}

In this work a Lagrangian descriptor approach to computing invariant manifolds 
is applied to an analysis of the roaming mechanism for chemical reaction dynamics. 
In the context of the CH$_4^+$ model of Chesnavich this means we want to 
uncover the phase space mechanism whereby the  hydrogen atom separates from the CH$_3^+$ core 
and later returns to the core before dissociating. 
We now review the dynamical definition of roaming that was introduced in \cite{mauguiere2014roaming}, 
which is based on periodic orbits as the phase space structures governing the dynamics of roaming.

\begin{figure}[ht]
 \centering
 \includegraphics[width=.5\textwidth]{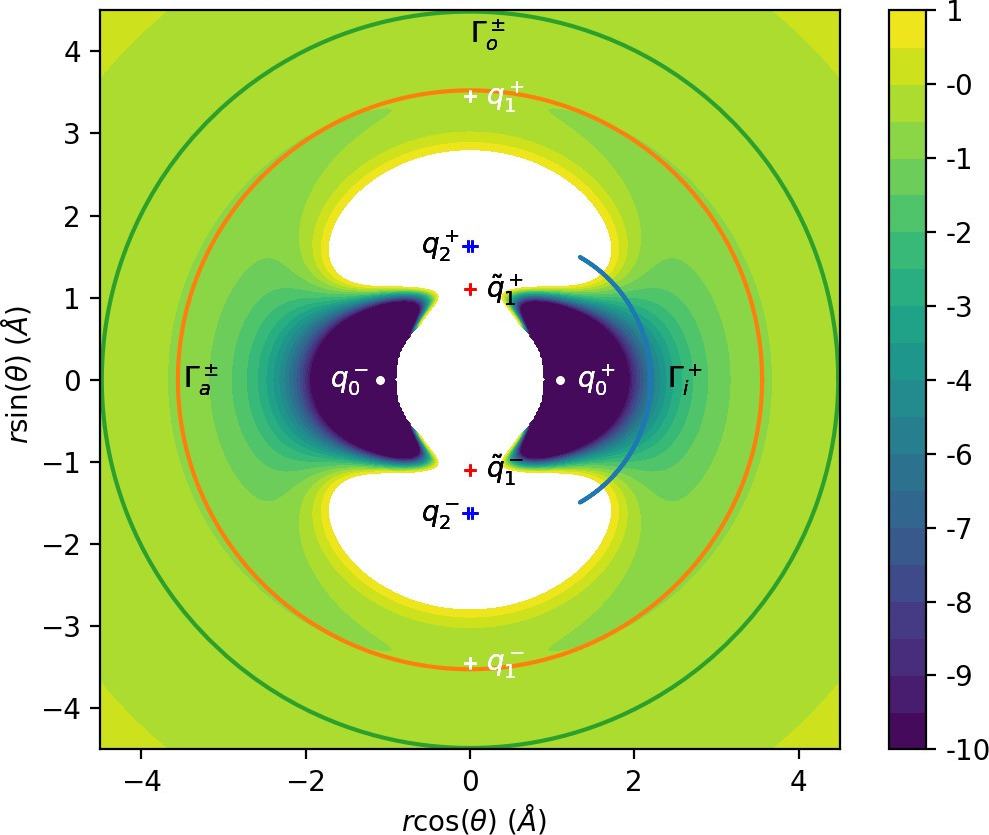}
 \caption{Configuration space projections of the inner ($\Gamma^i$), middle ($\Gamma^a$) and outer ($\Gamma^o$) periodic orbits for $E=5$. This figure is from \cite{krajnak2019isokinetic}.}
 \label{fig:orbitse5}
\end{figure}

In the relevant energy interval, $0\leq E\leq 5$, there are three  families of periodic orbits
that organize the roaming dynamics in phase space \cite{mauguiere2014multiple}. 
These periodic orbits are pairwise related by symmetry \eqref{eq:sym}, as 
shown in Fig. \ref{fig:orbitse5}. We will refer to these families as the inner ($\Gamma^i$), 
middle ($\Gamma^a$) and outer ($\Gamma^o$) periodic orbits. 
We refer to a continuum of periodic orbits parametrised by energy as a family of periodic orbits. 
We highlight the fact  that none of the orbits is directly related to a saddle point on the potential energy surface.

The significance of these periodic orbits for the roaming dynamics is as follows:

\begin{itemize}

 \item $\Gamma^i$: Delimits the potential wells that correspond to CH$_4^+$ isomers. 
 The orbits oscillate about the axes $\theta=0$ and $\theta=\pi$.
 
 \item $\Gamma^a$: Two rotational orbits with opposite orientation - one clockwise, one counter-clockwise - 
 that are crucial for defining roaming.
 
 \item $\Gamma^o$: Centrifugal barrier delimiting the region of dissociated states. 
 Two rotational orbits with opposite orientation - one clockwise, one counter-clockwise.
 
\end{itemize}

It was shown in \cite{krajnak2018phase} that roaming does not occur for $E\geq 2.5$. 
In the energy interval $0< E< 2.5$, all of the above-mentioned periodic orbits are unstable. 
We use these periodic orbits to define dividing surfaces. Denote by DS$^i$, DS$^a$ and DS$^o$ the 
set of all points $(r,\theta,p_r, p_\theta)$ on the energy surface 
$H(r,\theta,p_r, p_\theta)=E,$ whose configuration space projections 
$(r,\theta)$ coincide with the configuration space projections of 
$\Gamma^i$, $\Gamma^a$ and $\Gamma^o$ respectively. Due to the instability of the orbits, 
the resulting dividing surfaces satisfy the non-recrossing properties.

A roaming trajectory is then defined as a trajectory that crosses DS$^a$ an odd number of times between 
leaving the potential well and dissociating. An isomerising trajectory leaves the potential well, 
crosses DS$^a$ an even number of times and returns to either of the potential wells. 
A nonreactive trajectory originates in dissociated states and returns there after an even number of crossings of DS$^a$.

As explained in \cite{krajnak2018phase}, 
DS$^i$ consists of two spheres because $\Gamma^i$ is comprised of self-retracing (brake) orbits, 
while the rotational orbits $\Gamma^a$ and $\Gamma^o$ imply that DS$^a$ and DS$^o$ are tori \cite{mauguiere2016ozone}.
Each sphere can be divided using the corresponding periodic orbit into two hemispheres and 
each torus can be divided using both corresponding periodic orbits into two annuli \cite{mauguiere2016ozone}. 
All hemispheres and annuli are surfaces of unidirectional flux.  This implies that all 
trajectories leaving the potential well cross the same (outward) hemisphere of DS$^i$, 
while all trajectories entering the potential well cross the other (inward) hemisphere of DS$^i$.

Roaming can be described  as a transport problem in phase space. Every trajectory leaving the well 
must cross the outward hemisphere of DS$^i$ and every trajectory that dissociates must cross 
the outward annulus of DS$^o$. Dissociation of CH$_4^+$ is therefore 
equivalent to the transport of trajectories from the outward hemisphere of DS$^i$ to 
the outward annulus of DS$^o$. Roaming involves crossing the inward annulus of DS$^a$, 
since between two crossings of the outward annulus trajectories must cross the inward annulus and vice versa.

\begin{figure}[ht]
 \centering
 \includegraphics{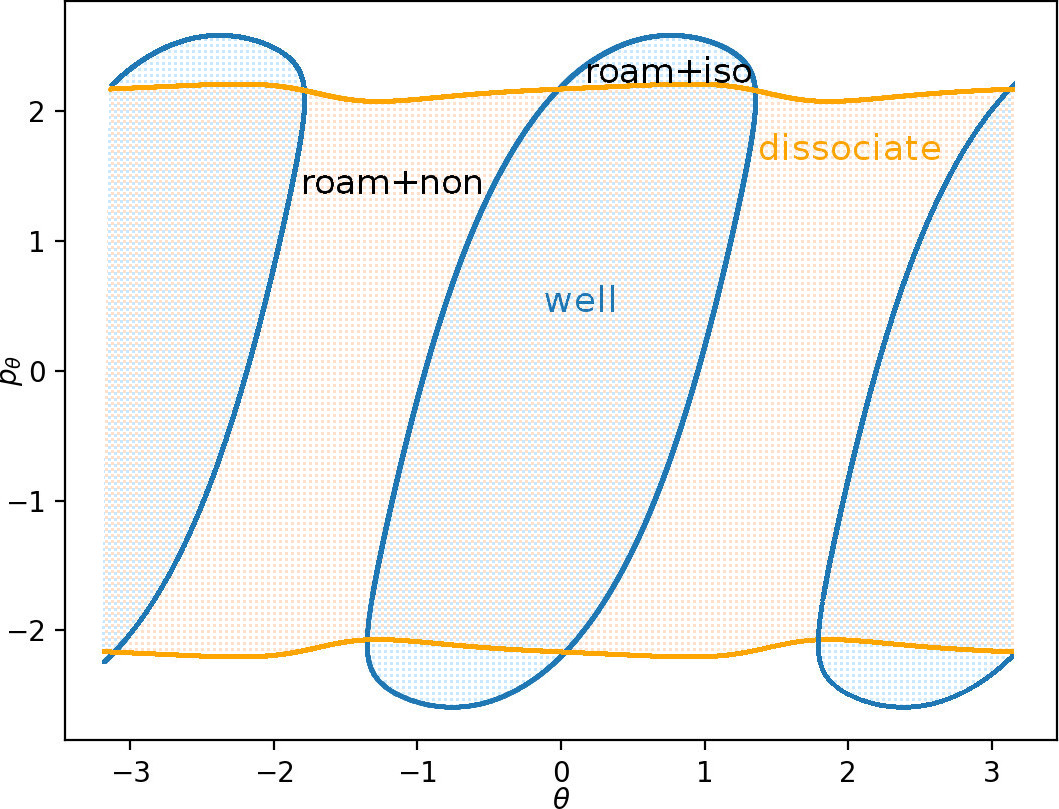}
 \caption{Intersection of invariant manifolds of $\Gamma^i$ (blue) and $\Gamma^o$ (orange) with the outward annulus of DS$^a$ for $E=1$. Trajectories that just left the potential well are shown in blue, immediately dissociating trajectories in orange. Roaming and isomerising trajectories in the blue area just left the well and do not dissociate immediately, while roaming and dissociating trajectories in the orange are dissociate immediately but did not just leave the potential wells.  This figure is from \cite{krajnak2019isokinetic}.}
 \label{fig:DSa}
\end{figure}

Transport of trajectories in the neighbourhood of an unstable periodic orbit (or NHIM in general) is 
governed by invariant manifolds of the orbit \cite{wwju, uzer2002geometry, WaalkensSchubertWiggins08, wiggins2016}. 
It was shown \cite{krajnak2018phase, krajnak2018influence} that the roaming phenomenon 
involves a heteroclinic intersection of the invariant manifolds of $\Gamma^i$ and $\Gamma^o$. 
The condition  $H(r,\theta=0,p_r>0, p_\theta=0)$ defines an invariant subsystem that consists of 
precisely one dissociating trajectory for every fixed $E>0$. Therefore if the invariant manifolds of 
$\Gamma^i$ and $\Gamma^o$ do not intersect, the former are contained in the interior of the latter 
and each trajectory leaving the potential well dissociates directly. 
An intersection assures that some trajectories leaving the well do not dissociate directly 
but return to DS$^a$ as illustrated in Fig. \ref{fig:DSa}. 
This is the phase space mechanism  for roaming and isomerisation in the Hamiltonian case.

\section{The Hamiltonian Isokinetic Thermostat for the Chesnavich Model}
\label{sec: HIK}

Here we apply the Hamiltonian isokinetic thermostat to the Chesnavich model. 
Details of the theory behind the Hamiltonian isokinetic thermostate model in this 
context can be found in \cite{krajnak2019isokinetic}. Here we just give the relevant equations with a brief description.

The isokinetic Hamiltonian $\cK$ for Chesnavich's model is
\begin{equation}
 \cK(r,\pi_r,\theta,\pi_\theta) = \frac{1}{2}e^{U}\left(\frac{\pi_r^2}{\mu}+\pi_\theta^2\left(\frac{1}{\mu r^2}+\frac{1}{I_{CH_3}}\right)\right) - \frac{1}{2}e^{-U},
\label{eq:chesKHam}
\end{equation}
where $U=U(r,\theta)$ is Chesnavich's potential energy and 
\begin{equation}
 \pi_r=e^{-U}p_r,\quad \pi_\theta=e^{-U}p_\theta.
 \label{eq:isotransf}
\end{equation}

The level set $\cK=0$ corresponds to the surface of constant kinetic energy
\begin{equation}
 \frac{1}{2}\frac{p_r^2}{\mu}+\frac{1}{2}p_\theta^2\left(\frac{1}{\mu r^2}+\frac{1}{I_{CH_3}}\right)=\frac{1}{2},
 \label{eq:constantT}
\end{equation}
in system \eqref{eq:chesHam}. We fix the kinetic energy at $\frac{1}{2}$, 
because according to \cite{krajnak2019isokinetic} the dynamics generated 
by \eqref{eq:chesKHam} is equivalent regardless of the value of kinetic energy
(up to a time rescaling). 
This statement remains true for other systems provided that the kinetic energy is a quadratic form in the momenta.

The equations of motion in the isokinetic system in terms of variables $(r, \pi_r, \theta, \pi_\theta)$ 
are given in \cite{krajnak2019isokinetic}.
To achieve greater numerical precision, it is preferable to integrate the equations of 
motion in the $(r,p_r,\theta,p_\theta)$ coordinates instead of the $(r,\pi_r,\theta,\pi_\theta)$
\cite{krajnak2019isokinetic}:
\begin{equation}
 \begin{split}
  \dot{r}&=\frac{p_r}{\mu},\\
 \dot{p}_r&= p_r\left(\frac{\partial U}{\partial r}\dot{r}+\frac{\partial U}{\partial \theta}\dot{\theta}\right) + \frac{1}{\mu r^3}p_\theta^2 -\frac{\partial U}{\partial r},\\
 \dot{\theta}&=p_\theta\left(\frac{1}{\mu r^2}+\frac{1}{I_{CH_3}}\right),\\
 \dot{p}_\theta&= p_\theta\left(\frac{\partial U}{\partial r}\dot{r}+\frac{\partial U}{\partial \theta}\dot{\theta}\right) -\frac{\partial U}{\partial \theta}, 
 \label{eq:chesKHameqreduced}
 \end{split}
\end{equation}

\noindent
where we used the isokinetic constraint \eqref{eq:constantT} equivalent to $\cK=0$.

\begin{figure}[ht]
 \centering
 \includegraphics[width=.6\textwidth]{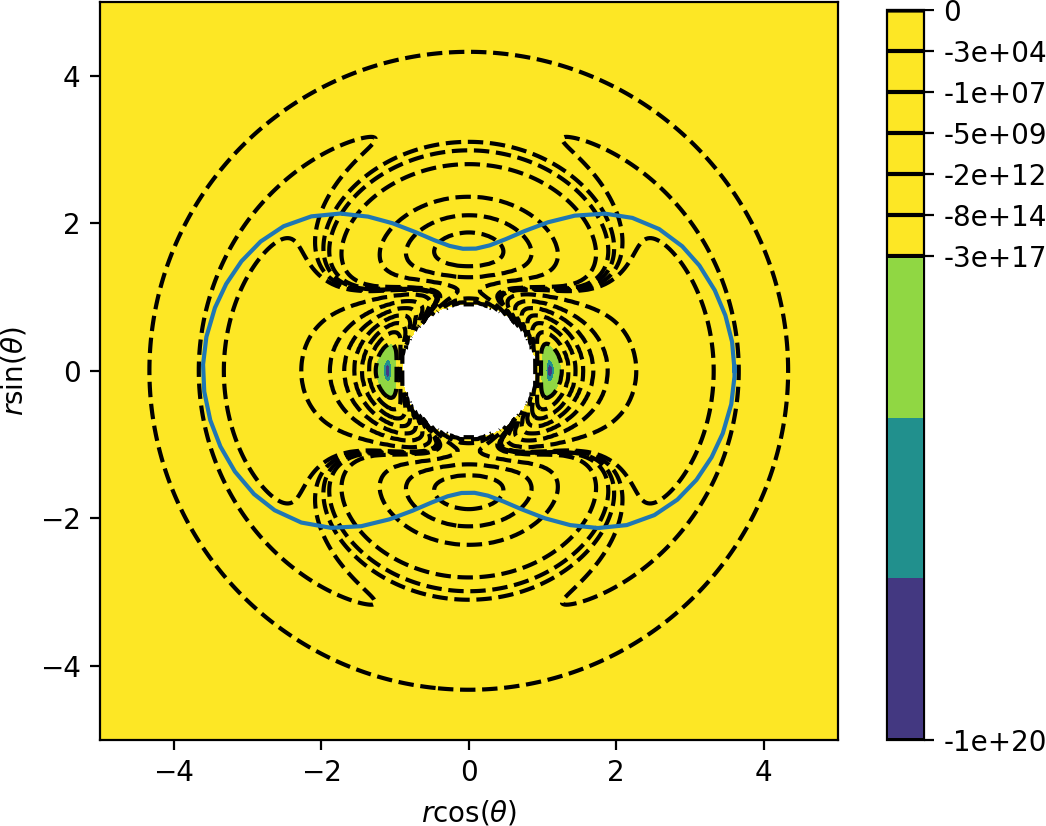}
 \caption{Inner periodic orbit on the potential energy surface $- \frac{1}{2}e^{-U}$.  
 This figure is from \cite{krajnak2019isokinetic}.}
 \label{fig:disI}
\end{figure}

The potential $- \frac{1}{2}e^{-U}$ has the same critical points and characteristics as $U$, 
but the wells are considerably deeper and have steeper walls. In contrast to the microcanonical case, 
the isokinetic model only possesses two periodic orbits with period $2\pi$ and due 
to constant nonzero kinetic energy does not admit self-retracing orbits (also referred to as brake orbits)
such as $\Gamma^i$ introduced in Section \ref{sec:roamingChes}. 
One of the periodic orbits delimits the potential wells, see Fig. \ref{fig:disI}; we therefore refer to it as the inner orbit.

The outer orbit, beyond which trajectories do not return and corresponds to the dissociated state of the molecule, 
is associated with a centrifugal barrier. It is nearly rotationally symmetric and in double precision has a constant radius $r=13.4309241401910709$. 

Its existence can be proven using a similar argument as in the original system \cite{krajnak2018phase}: 
suppose $r$ is sufficiently large so that $U$ is effectively independent of $\theta$. Denote $r_{po}$ the solution of
\begin{equation}
\frac{1}{\mu r_{po}^3}p_\theta^2 -\frac{\partial U}{\partial r}=0.
\end{equation}
Then the equations \eqref{eq:chesKHameqreduced} admit a rotationally symmetric periodic orbit 
with $\dot{\theta}=const$, provided
\begin{equation}
 \begin{split}
  \dot{r}&=0,\\
 \dot{p}_r&=0.
 \end{split}
\end{equation}

This is satisfied by the initial condition $r=r_{po}$, $p_r=0$ and $p_\theta$ 
given implicitly by $\cK=0$ for any $\theta$. The existence of $r_{po}$ is guaranteed 
for the potential $U$ and any other potential with leading order term $-cr^{-(2+\varepsilon)}$ for 
large $r$, with $c>0$ and $\varepsilon>0$.

Both these orbits are unstable, with the largest eigenvalue of the inner orbit 
under the return map being of the order $10^{21}$. This large instability poses a 
serious challenge to calculation of its invariant manifolds that guide trajectories in phase space. 
We therefore employ the method of Lagrangian descriptors to find these invariant manifolds.


\section{Lagrangian descriptors and invariant manifolds}
\label{sec:LDIM}
\subsection{Introduction to Lagrangian descriptors}
Lagrangian descriptors are a very successful trajectory diagnostic for revealing phase space structures 
in dynamical systems. The method was originally developed for analyzing Lagrangian transport phenomena 
in fluid dynamics \cite{madrid2009ld}. but the utility and applicability of the method has recently 
been recognized in chemistry 
\cite{craven2016deconstructing, craven2017lagrangian, craven2015lagrangian, junginger2016transition, junginger2017chemical, feldmaier2017obtaining, junginger2016transition, feldmaier2017obtaining, patra2018detecting}. 
The method is simple to implement computationally, the interpretation 
in terms of trajectory behaviour is clear, and it provides a `high resolution' method for exploring 
high-dimensional phase space with low-dimensional slices \cite{ demian2017detection, naik2019a, naik2019b, garcia2019tilting}. 
It applies to Hamiltonian and non-Hamiltonian systems \cite{lopesino2017} 
and to systems with arbitrary, even stochastic, time-dependence \cite{balibrea2016lagrangian} 
Moreover, Lagrangian descriptors can be applied directly to data sets, without the need of 
an explicit dynamical system \cite{mendoza2014lagrangian}.

\subsection{Lagrangian descriptors minimised by invariant manifolds}
\label{sec:newLD}
Complementary to existing applications of Lagrangian descriptors as an exploratory method, 
we use Lagrangian descriptors to find invariant manifolds in a system 
where some invariant structures have already been identified. 
In the present case we use the parametrisation of a periodic orbit to construct a Lagrangian descriptor 
that will be locally minimised by its invariant manifolds. 

For a Hamiltonian system with $2$ degrees of freedom and $(q_1,p_1,q_2,p_2)$ as phase space coordinates, 
a $1$ dimensional periodic orbit $\Gamma$ can be viewed as an intersection of three hypersurfaces:
the surface of constant energy it lives on and two other distinct hypersurfaces that 
parametrise $\Gamma$, $f_1(q_1,p_1,q_2,p_2)=0$ and $f_2(q_1,p_1,q_2,p_2)=0$. 
For a system with zero total angular momentum, these equations can correspond to 
the configuration space projection and parametrisation of one of the momenta.

The integral
\begin{equation}
 \int\limits^{\infty}_{0} |\dot{f_i}(q_1,p_1,q_2,p_2)| dt,
 \label{eq:LDinfty}
\end{equation}
then vanishes along the periodic orbit for $i=1,2$.

As $t\rightarrow\infty$, $f_1,f_2\rightarrow 0$ along the periodic orbit stable invariant manifolds. 
Consequently of all points on the energy surface, \eqref{eq:LDinfty} is locally minimised by 
those that lie on the stable invariant manifold of $\Gamma$. Necessarily this statement remains true for 
the finite time integral
\begin{equation}
 \int\limits^{\tau}_{0} |\dot{f_i}(q,p)| dt,\qquad i=1,2,
\end{equation}
provided $\tau>0$ is large enough. We can use the same argument for the unstable invariant manifold and the Lagrangian descriptor
\begin{equation}
 \int\limits^{0}_{-\tau} |\dot{f_i}(q,p)| dt,\qquad i=1,2.
\end{equation}

Depending on the strength of instability of $\Gamma$, small values of $\tau$ may suffice to find 
an approximate location of the invariant manifolds and to establish whether or not they intersect, which is a crucial characteristic for the presence of roaming.

\subsection{Roaming in the Isokinetic Chesnavich Model}
In \cite{krajnak2019isokinetic}, we established the presence of roaming and 
classified dynamical behaviour using escape time analysis and Lagrangian descriptors. 
Classes of qualitatively different dynamical behaviour are separated by invariant manifolds, 
since these guide trajectories across phase space bottlenecks.

Due to the large instability of the inner orbit, it was not previously 
possible to compute its invariant manifolds to create an analogue of Fig. \ref{fig:DSa}. 
The Lagrangian descriptors defined in Sec. \ref{sec:newLD} uncover trajectories in the proximity 
of invariant manifolds and we thereby avoid having to compute invariant manifolds at all.

As for the microcanonical version of Chesnavich's model, we are interested in the following manifolds:
\begin{itemize}
 \item the unstable manifold of the inner orbit $W_i^u$, which guides trajectories out of the wells into the interaction region between DS$^i$ and DS$^o$,
 \item the stable manifold of the outer orbit $W_o^s$, which guides trajectories out of the interaction region towards dissociation.
\end{itemize}

We now construct the Lagrangian descriptors described in Sec. \ref{sec:newLD}.
The inner periodic orbit is parametrised by
\begin{equation}
f_1(r, \theta) = r-\bar{r}(\theta)=0,
\end{equation}
where
\begin{equation}
\bar{r}(\theta)=c_0+c_1\cos(2\theta)+c_2\cos(4\theta)+c_3\cos(6\theta)+c_4\cos(8\theta)+c_5\cos(10\theta),
\end{equation}
with constants
\begin{equation}
 \begin{split}
  c_0=2.78147867,\ c_1=0.98235111,\ c_2=-0.17161848,\\
  c_3=-0.00486657,\ c_4=0.01628185,\ c_5=-0.00393858,
 \end{split}
\end{equation}
and
\begin{equation}
f_2(\theta, p_r) = p_r-\bar{p}_r(\theta)=0,
\end{equation}
where $\bar{p}_r(\theta)$ is a function defined for $\theta\in (-\frac{\pi}{2},\frac{\pi}{2}]$ and periodically extended to $\theta\in\mathbb{R}$, given by
\begin{equation}
\bar{p}_r(\theta)=d_0\theta + d_1\theta^3 + d_2\theta^5 + d_3\theta^7 + d_4\theta^9 + d_5\theta^{11},
\end{equation}
with constants
\begin{equation}
 \begin{split}
  d_0=-1.06278495,\ d_1=-0.42089795,\ d_2=1.38849679,\\
  d_3= -1.11654771,\ d_4=0.40789372,\ d_5=-0.05122644.
 \end{split}
\end{equation}
Note that the periodic orbit need only be parametrised for a single value of $\cK$,
$\cK = 0$.

To identify $W_i^u$, we use the Lagrangian descriptor
\begin{equation}
 LD_i=\int\limits^{0}_{-\tau} |\dot{r}-\bar{r}^{'}(\theta)\dot{\theta}| dt.
 \label{eq:LDi}
\end{equation}
Note that we obtain identical results using
\begin{equation*}
 \int\limits^{0}_{-\tau} |\dot{p}_r-\bar{p}_r^{'}(\theta)\dot{\theta}| dt.
\end{equation*}

We approximate the outer periodic orbit by a constant radius $r=13.4309241401910709$ and 
its stable invariant manifold $W_o^s$ locally minimises
\begin{equation}
 LD_o=\int\limits^{\tau}_{0} |\dot{r}| dt.
 \label{eq:LDo}
\end{equation}

\begin{figure}[ht]
 \centering
 \includegraphics[width=0.49\textwidth]{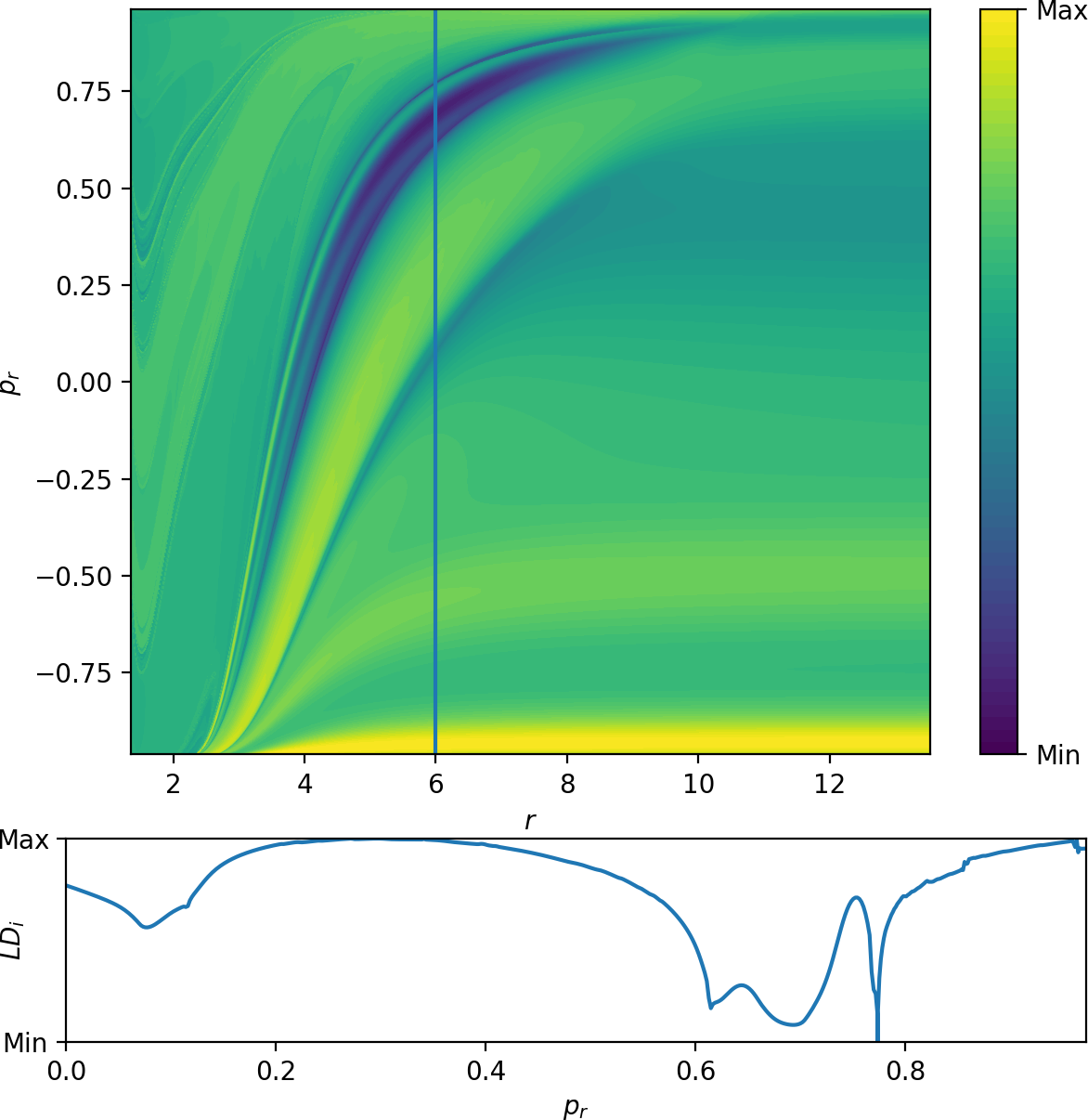}
 \includegraphics[width=0.49\textwidth]{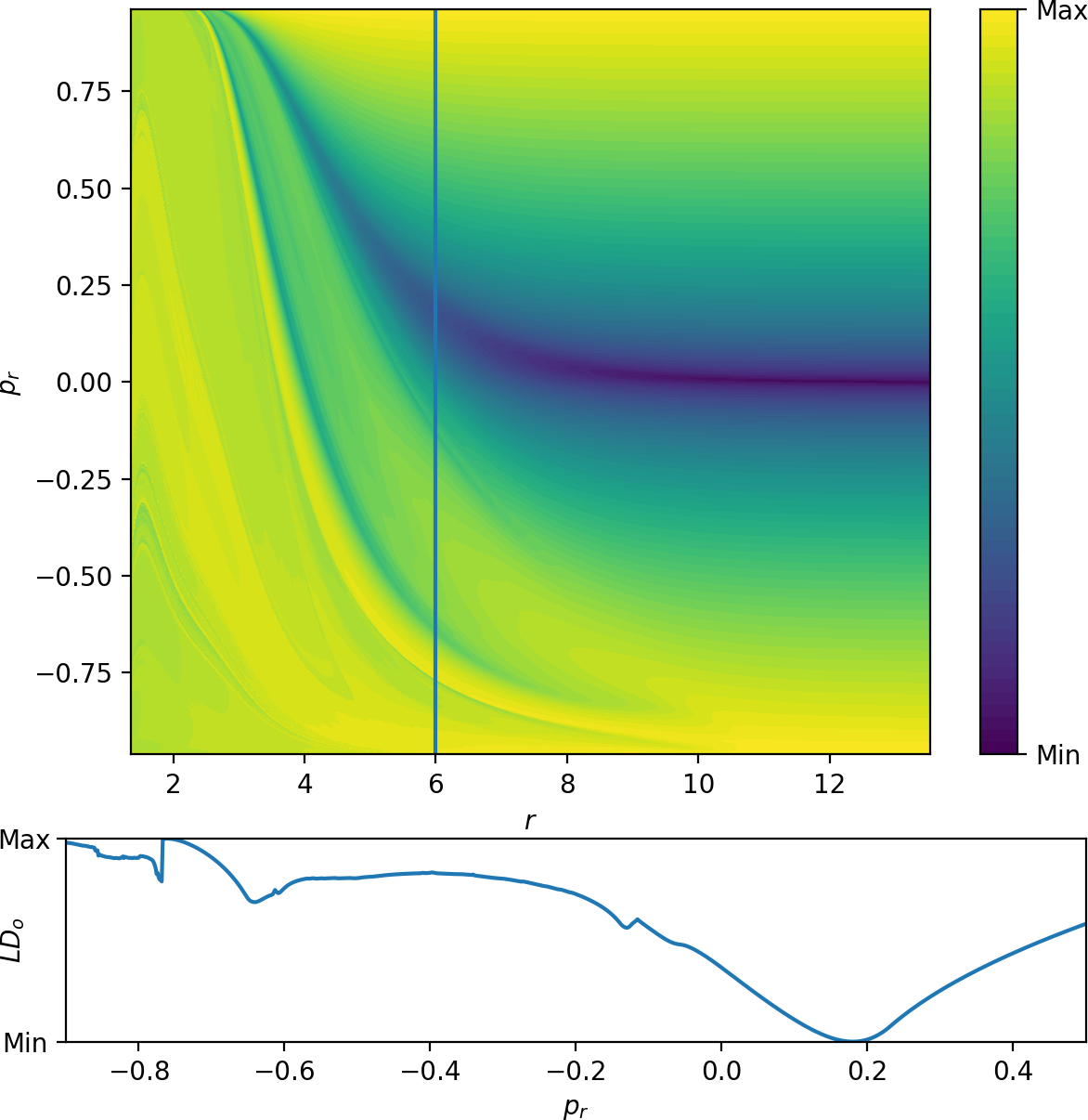}
 \caption{Lagrangian descriptors $LD_i$ \eqref{eq:LDi} for $\tau=8$ (left) and $LD_o$ \eqref{eq:LDo} for $\tau=8$ (right) for initial conditions on $\theta=0$, $\dot{\theta}>0$ and their profiles for $r=6$.}
 \label{fig:manifs_th0} 
\end{figure}

In \cite{krajnak2019isokinetic}, we found that singularities in the escape time plot
matched significant changes in Lagrangian descriptor plots. Both these features indicated 
the presence of invariant manifolds. In Fig. \ref{fig:manifs_th0} we show 
$LD_i$ \eqref{eq:LDi} and $LD_o$ \eqref{eq:LDo} on the surface $\theta=0$, $\dot{\theta}>0$, 
where the latter is similar to the plots in \cite{krajnak2019isokinetic} and the former shows similar features 
reflected about $p_r=0$ due to the opposite time direction in the definition of $LD_i$. 
Both plots show LD values on a uniform $400\times400$ grid. A higher density 
improves accuracy but does not yield additional qualitative information. 
The accompanying sections for $r=6$ illustrate the nature of the local minima.

The presence of roaming is only possible if $W_i^u$ and $W_o^s$ intersect \cite{krajnak2018phase}. 
Fig. \ref{fig:manifs_th0} provides direct evidence of the intersection - the two panels 
show $W_i^u$ and $W_o^s$ using different Lagrangian descriptors on the same surface of section. 
These manifolds necessarily intersect.

\begin{figure}[ht]
 \centering
 \includegraphics[width=0.49\textwidth]{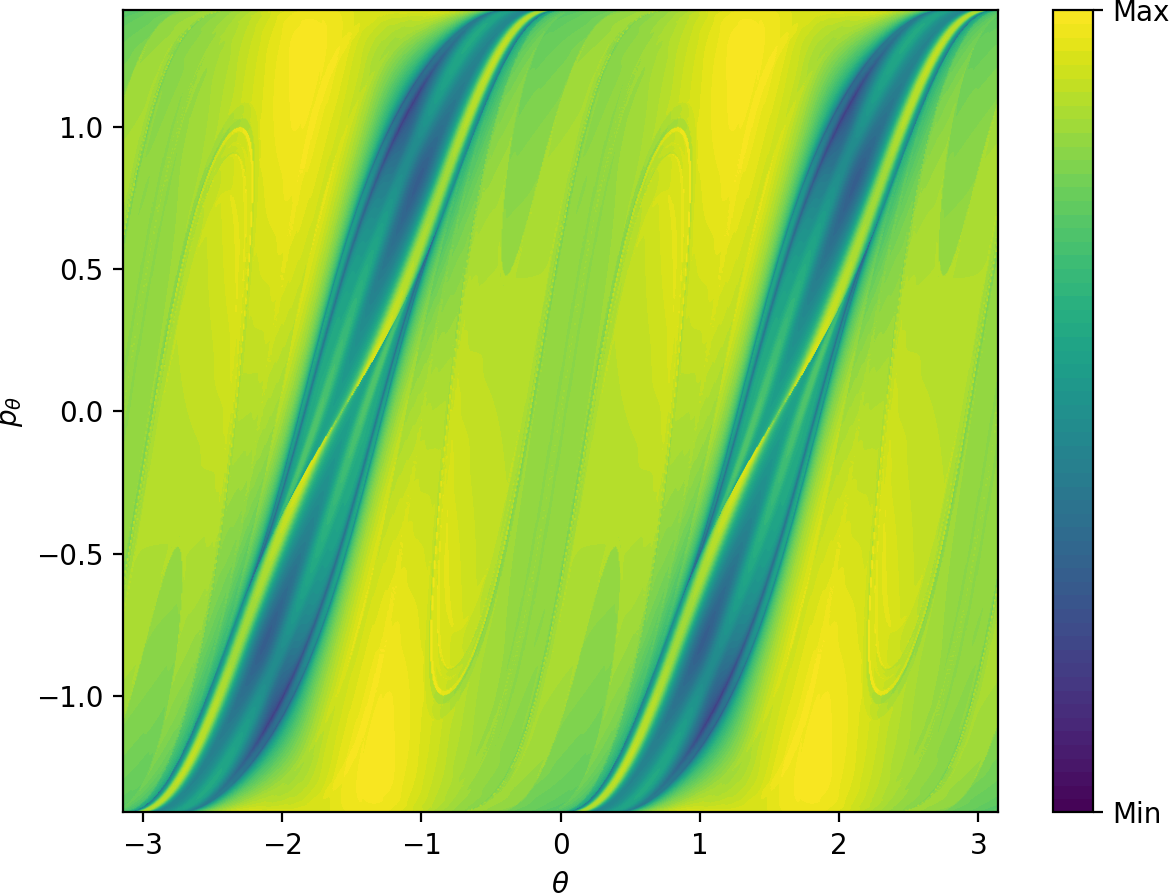}
 \includegraphics[width=0.49\textwidth]{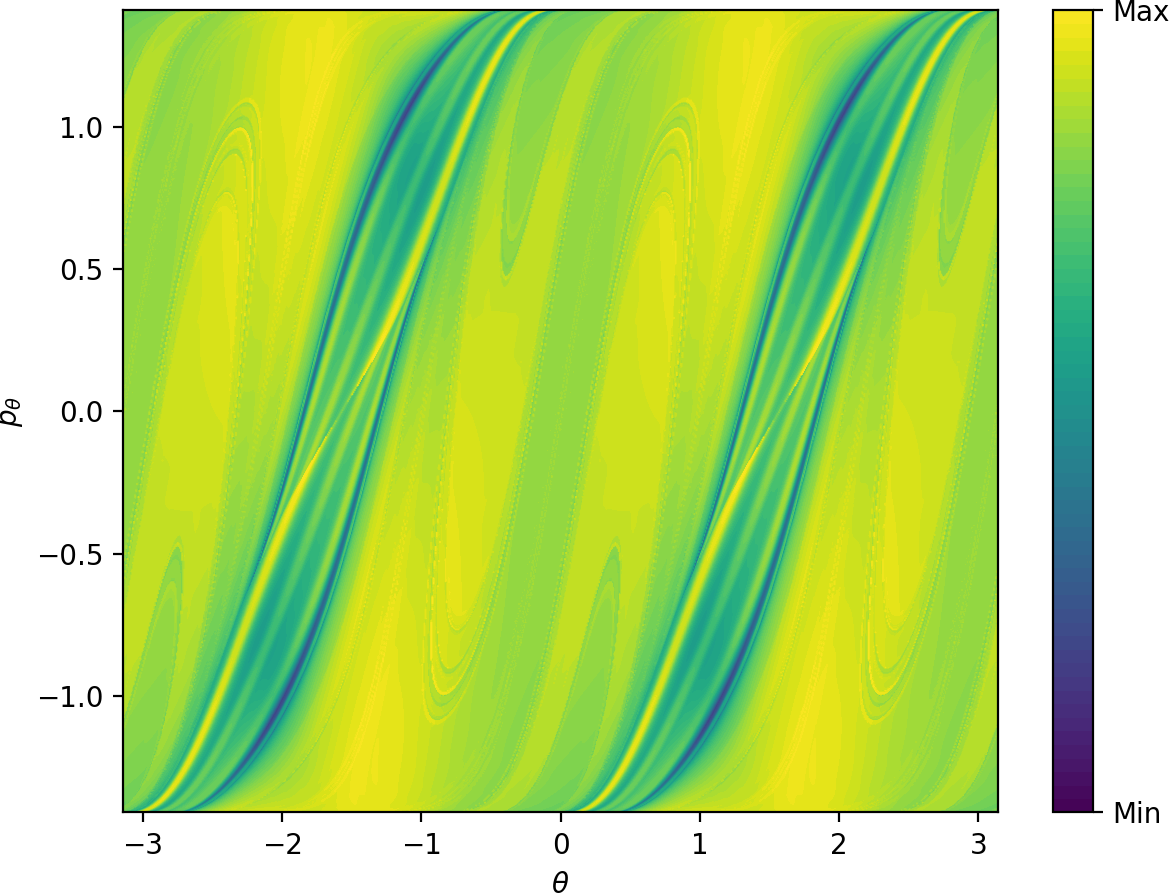}
 \caption{Lagrangian descriptor $LD_i$ \eqref{eq:LDi} for $\tau=6$ (left) and $\tau=8$ (right) for initial conditions on $r=3.6$, $\dot{r}>0$.}
 \label{fig:manifsI_r36} 
\end{figure}

Further information on the geometry of the intersection of $W_i^u$ and $W_o^s$ can be 
seen using a surface of section analogous to the outward annulus of DS$^a$ used in the Hamiltonian case.

\begin{figure}[ht]
 \centering
 \includegraphics[width=0.49\textwidth]{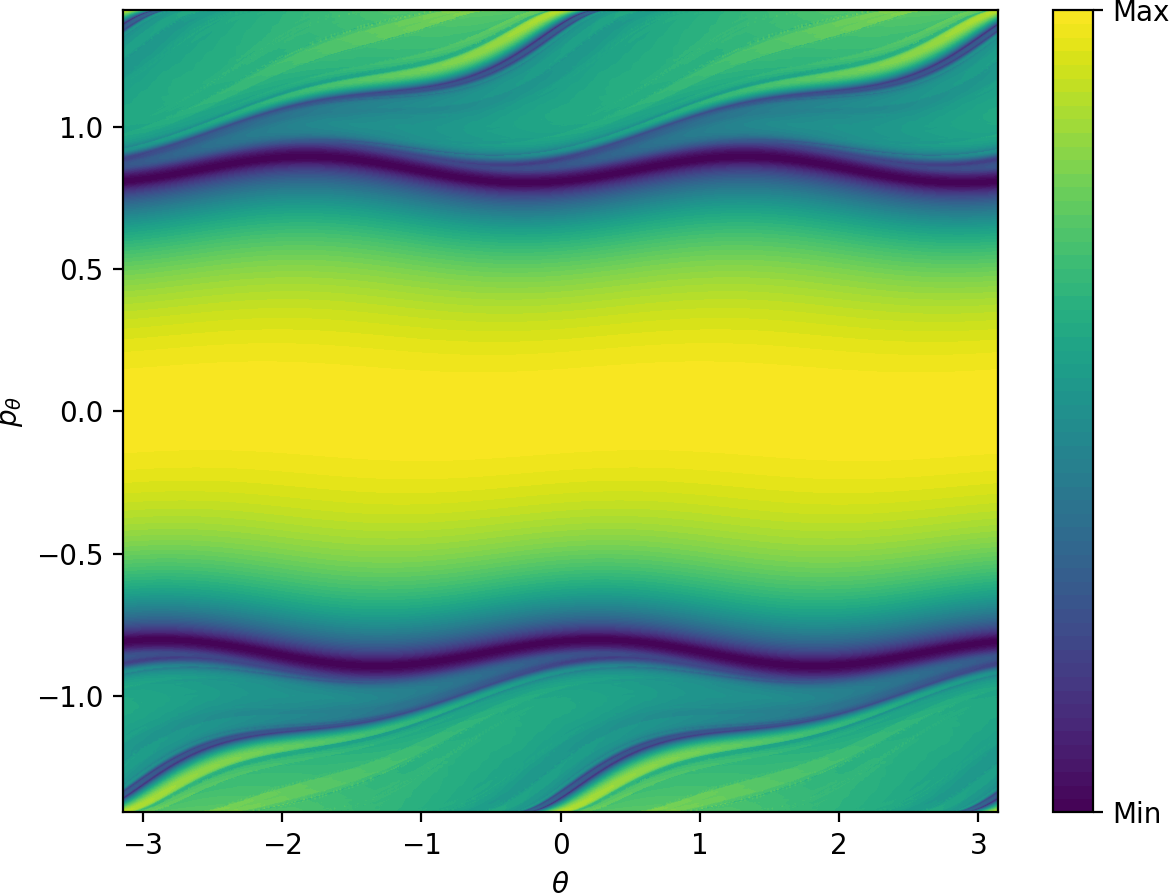}
 \includegraphics[width=0.49\textwidth]{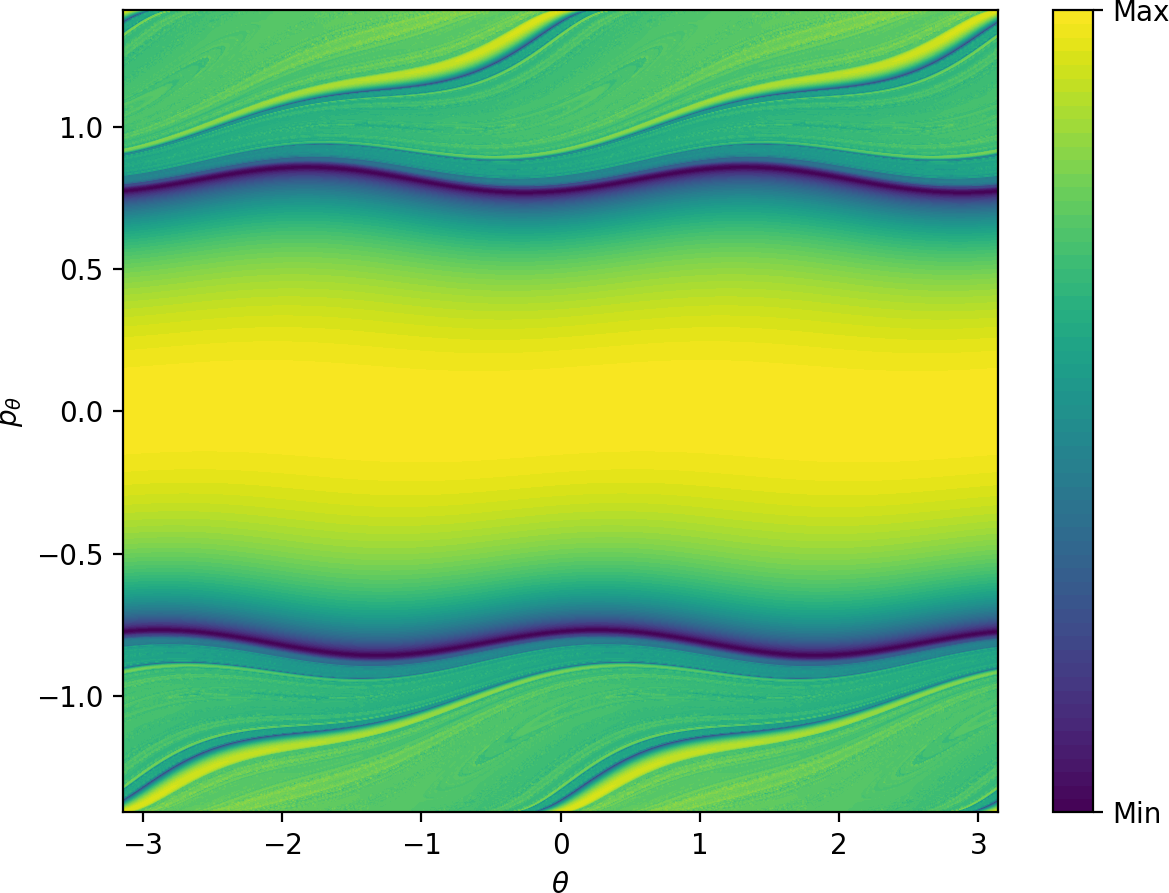}
 \caption{Lagrangian descriptor $LD_o$ \eqref{eq:LDo} for $\tau=10$ (left) and for $\tau=20$ (right) for initial conditions on $r=3.6$, $\dot{r}>0$.}
 \label{fig:manifsO_r36}
\end{figure}

Due to the absence of other $2\pi$-periodic orbits, one cannot expect to easily obtain 
a surface of section that is transversal to the flow. For simplicity, we use 
a surface of section defined by the condition $r=3.6$, $\dot{r}>0$. 
Figures \ref{fig:manifsI_r36} and \ref{fig:manifsO_r36} 
show $LD_i$ \eqref{eq:LDi} and $LD_o$ \eqref{eq:LDo} on this surface. 
Areas of low values of $LD_i$ and $LD_o$ correspond to trajectories that are
asymptotic to the respective periodic orbits. Note that the structures are well 
defined visually for times that are of the same order of magnitude as the periodic orbits 
themselves - periods of the inner and outer orbits are $11.84$ and $9.61$ respectively.

Similarly to the microcanonical case, the toriodal cylinder $W_o^s$ intersects $r=3.6$, $\dot{r}>0$ 
transversally in two circles, one with $p_\theta>0$, one with $p_\theta<0$.

$W_i^u$ has the same geometry as $W_o^s$, but it does not intersect the surface of section transversally. 
This can be seen by the manifold attaining maximal $p_\theta$ permitted by kinetic energy $\frac{1}{2}$, when 
$p_r$ vanishes and $\dot{r}>0$ is violated. As a result, instead of two circles, the two unstable invariant 
manifolds intersect the surface of section in four S-shapes. 
The initial condition $\theta=p_\theta=0$ leads to radial dissociation, 
and it is surrounded by trajectories that originate in the well. This area is bounded by two of 
the S-shapes - the S-shape passing closest to $\theta=1$, $p_\theta=0$ corresponds 
to the inner periodic orbit with $p_\theta>0$. The S-shape closest to $\theta=-1$, $p_\theta=0$ 
corresponds to the inner periodic orbit with $p_\theta<0$. The other S-shapes are related by symmetry.

\begin{figure}[ht]
 \centering
 \includegraphics{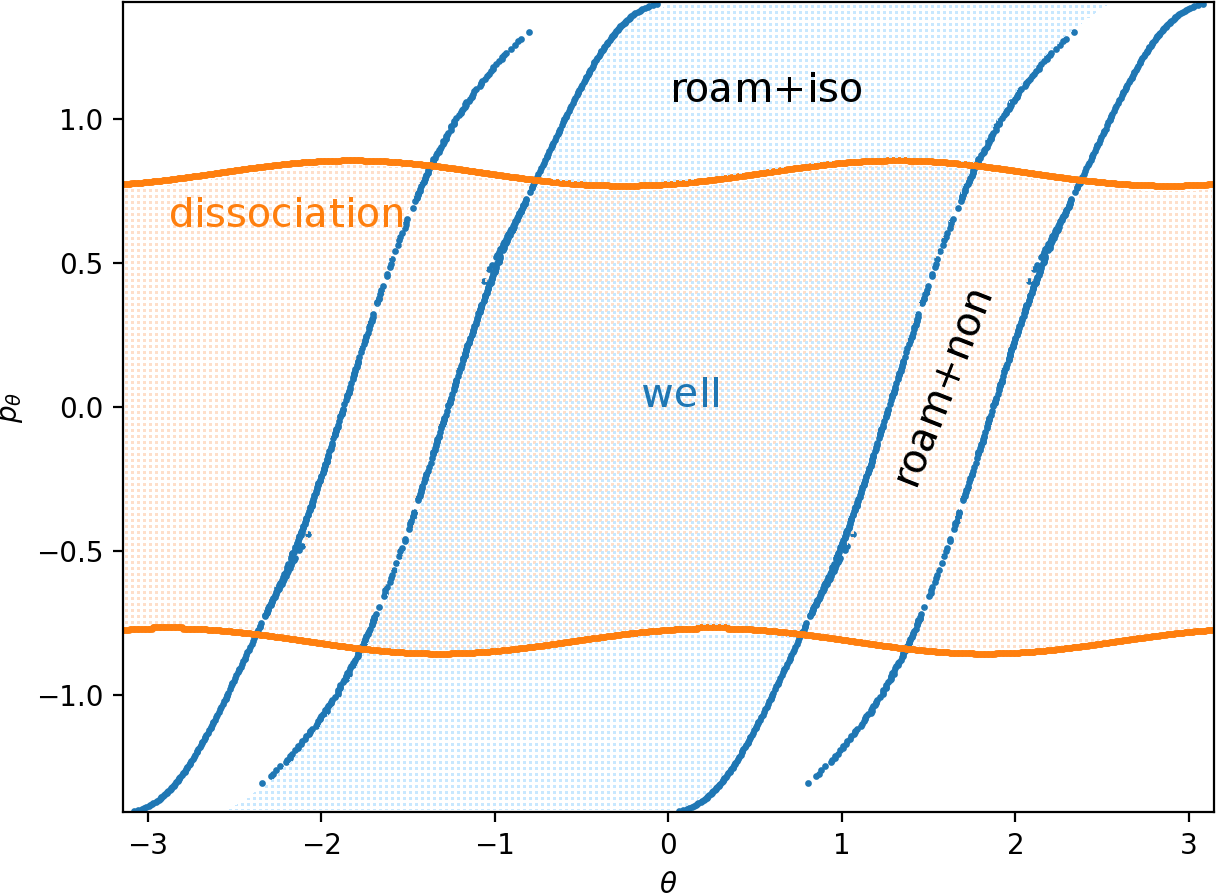}
 \caption{Approximations of invariant manifolds on $r=3.6$, $\dot{r}>0$ extracted from $LD_i$ for $\tau=6$ shown in Fig. \ref{fig:manifsI_r36} and $LD_o$ for $\tau=20$ shown in Fig \ref{fig:manifsO_r36} that separate different classes of dynamics. For details see text.}
 \label{fig:intersection}
\end{figure}

Fig. \ref{fig:intersection} shows superposed manifolds approximated from $LD_i$ ($\tau=6$) and $LD_o$ ($\tau=20$) values. 
Minima of $LD_o$ for fixed values of $\theta$ reveal $W_o^s$ immediately.

Due to the high instability of the inner orbit, locating $W_i^u$ using minima of $LD_i$ requires 
a very dense sampling with a correspondingly high computational cost. 
(Note how narrow the interval for finding the minimum on the profile of $LD_i$ on $\theta=0$, $r=6$ 
in Fig. \ref{fig:manifs_th0} is.)
Instead we can take advantage of the steepnes of the gradient near the invariant manifold and consider 
maxima of $\Delta_\theta LD_i$ on a $400\times400$ grid for fixed $p_\theta$ values. 
This way we recover the four S-shapes due to $W_i^u$.
Applying a cutoff at high values of $LD_i$ before calculating $\Delta_\theta LD_i$ delivers a cleaner image;
in Fig. \ref{fig:intersection} we used a cutoff at $\frac{1}{2}(\text{max}_{LD_i}+\text{min}_{LD_i})$.

Aside from $W_i^u$, Fig. \ref{fig:intersection} shows traces of two nearly linear 
segments near $W_i^u$ pointing towards the points $\theta=\frac{\pi}{2}$, 
$p_\theta=0$ and $\theta=\frac{3\pi}{2}$, $p_\theta=0$ on the surface of section. 
These are initial conditions of
trajectories asymptotic to the manifolds $\theta=\frac{\pi}{2}$, $p_\theta=0$, $p_r>0$ and $\theta=\frac{3\pi}{2}$, 
$p_\theta=0$, $p_r>0$. (Due to conservation of kinetic energy, these manifolds are 
distinct from $\theta=\frac{\pi}{2}$, $p_\theta=0$, $p_r<0$ and $\theta=\frac{3\pi}{2}$, $p_\theta=0$, $p_r<0$.)

Having identified $W_i^u$ and $W_o^s$, we can observe 
that these invariant manifolds guide trajectories in the same way as in the microcanonical model. 
Fig. \ref{fig:intersection} is analogous to Fig. \ref{fig:DSa}. 
The area between the two circles of $W_o^s$ is crossed by trajectories that are led by $W_o^s$ to dissociation. 
The area containing $\theta=p_\theta=0$ between two of the S-shapes of $W_i^u$ is crossed by 
trajectories that originate in the well.

Clearly the intersection of these two areas contains directly dissociating initial conditions. 
The complement between two of the S-shapes contains isomerising and roaming trajectories; 
in fact it is here the trajectories cross the surface of section for the first time after leaving the well. 
Similarly we observe a band between the two circles of $W_o^s$ which is crossed by trajectories that do 
not just leave the well. These may be roaming trajectories that originate in the well 
and recross the surface or nonreactive trajectories that never enter the well at all. 
This band marks the last crossing of the surface by these trajectories before dissociation. 
The remaining area corresponds to trajectories that stay in the interaction region before entering 
one of the wells or dissociating.

It is important to note that all classes of dynamics are separated by an intricate fractal 
structure made up of invariant manifolds, as shown in \cite{krajnak2019isokinetic}. 
The structure becomes visible when integrating any Lagrangian descriptor over a long time interval, 
see the details of Fig. \ref{fig:manifsI_r36} and \ref{fig:manifsO_r36}.

\section{Conclusions and Outlook}
\label{sec:conc}
We have shown how to compute stable and unstable invariant manifolds for which traditional methods fail due to high instability of the corresponding periodic orbit. 
For this purpose we propose the construction of a Lagrangian descriptor defined using 
an explicit parametrisation of the periodic orbit. 
The method is simple and easy to implement because it does not 
require calculating eigendirections of the orbit.
Additional advantages follow from lower demands on accuracy.

Invariant manifolds obtained this way enabled us to compare roaming in Chesnavich's CH$_4^+$ model and 
the isokinetic version of Chesnavich's CH$_4^+$ model subject to a Hamiltonian isokinetic theromostat. 
We conclude that regardless of detailed differences in phase space structures, 
the invariant structures responsibles for roaming bear remarkable similarities despite the 
nonholonomic constraint of constant kinetic energy.

\section*{Acknowledgments}
\noindent We acknowledge the support of  EPSRC Grant no. EP/P021123/1 and Office of Naval Research (Grant No.~N000141712220).

\newcommand{\etalchar}[1]{$^{#1}$}
\def\cprime{$'$}

\end{document}